        \documentstyle{article}
        \newcommand{\al}{\alpha}

        \newcommand{\lam}{\lambda}

        \newcommand{\om}{\omega}

        \newcommand{\Sig}{{\mathit \Sigma}}

        \font\smbbb=msbm5 
        \font\bbb=msbm7 
        \font\BBB=msbm10 
        
        \newcommand{\ZZ}{{\mbox{\BBB{Z}}}}

        \newcommand{\re}{{\mbox{\bbb{R}}}}
        \newcommand{\RE}{{\mbox{\BBB{R}}}}
        \newcommand{\PP}{{\mbox{\BBB{P}}}}
        \newcommand{\smco}{{\mbox{\smbbb{C}}}}
        \newcommand{\co}{{\mbox{\bbb{C}}}}
        \newcommand{\CO}{{\mbox{\BBB{C}}}}
        \newcommand{\CP}{\PP}
        \newcommand{\CM}{\CO^{\times}}
        \newcommand{\cm}{{\co^{\times}}}
        \font\frak=eufm10 at 11 pt
        \newcommand{\g}{\mbox{\frak{g}}}

        \newcommand{\gt}{\mbox{\frak{t}}}

        \newcommand{\OO}{{\cal O}}

        \newcounter{sect}
        \newcounter{subsect}
        \newcommand{\sect}[1]{\vspace{2ex}
                \addtocounter{sect}{1}\setcounter{subsect}{0}
                \begin{flushleft}
                {{\large\bf \arabic{sect}. {#1}}}
                \end{flushleft}
                \setcounter{thm}{0}
                \setcounter{equation}{0}
                \def\theequation{\arabic{sect}.\arabic{equation}}
                \def\thefigure{\arabic{sect}.\arabic{figure}}}
        
        \newtheorem{thm}{Theorem}[sect]

        \newtheorem{cor}[thm]{Corollary}
        \newtheorem{assump}[thm]{Assumption}
        
        \newtheorem{ex}[thm]{Example}
        \newtheorem{rmk}[thm]{Remark}
        \newcommand{\proof}[1]{\noindent {\em Proof.}$\quad$ {#1} $\hfill\Box$}
        \newcommand{\be}{\begin{equation}}
        \newcommand{\ee}{\end{equation}}
        \newcommand{\bea}{\begin{eqnarray}}
        \newcommand{\eea}{\end{eqnarray}}
        \newcommand{\nno}{\nonumber \\}
        \newcommand{\sep}[1]{\!\!\!\! &{#1}& \!\!\!\! }
        \newcommand{\eq}{\sep{=}}
        \newcommand{\vc}{\sep{ }}

        \newcommand{\bra}{\langle}
        \newcommand{\ket}{\rangle}

        \newcommand{\medwedge}{\mbox{\fontsize{12pt}{0pt}\selectfont $\wedge$}}

        \newcommand{\e}[1]{e^{{#1}}}
        \newcommand{\ii}{\sqrt{-1}}
        \newcommand{\ch}{{\,\mathrm{char}\,}}

        \newcommand{\set}[2]{\{{#1}\,|\,{#2}\}}
        
        \newcommand{\zero}{{\{0\}}}

        \newcommand{\pq}{^{pq}}

        \newcommand{\ka}{K\"ahler }
        \newcommand{\bb}{Bia\l ynicki-Birula }

        \newcommand{\two}[4]{\left\{    \begin{array}{ll}
                                        {#1}, & {\mbox{if }} {#2}, \\
                                        {#3}, & {\mbox{if }} {#4}
                                        \end{array}     \right.}

\begin{document}
$\!\,{}$

        \vspace{-5ex}

        \begin{flushright}
{\tt math.SG/9809192} (September, 1998)\\
Adelaide IGA preprint
        \end{flushright}

\vspace{5ex}

        \begin{center}
{\Large\bf A Note on Higher Cohomology Groups of K\"ahler Quotients}\\

        \vspace{4ex}
        {\large\rm Siye Wu}\footnote{E-mail address:
        {\tt swu@maths.adelaide.edu.au}}

        {\em Department of Pure Mathematics, University of Adelaide, 
        Adelaide, SA 5005, Australia}

        \end{center}

        \vspace{3ex}

        \begin{quote}
{\small {\bf Abstract.}
Consider a holomorphic torus action on a possibly non-compact K\"ahler
manifold.
We show that the higher cohomology groups appearing in the geometric
quantization of the symplectic quotient are isomorphic to the invariant
parts of the corresponding cohomology groups of the original manifold.
For non-Abelian group actions on compact K\"ahler manifolds,
this result was proved recently by Teleman and by Braverman.
Our approach is applying the holomorphic instanton complex
to the prequantum line bundles over the symplectic cuts.
We also settle a conjecture of Zhang and the present author on the exact 
sequence of higher cohomology groups in the context of symplectic cutting.}
        \end{quote}

        \vspace{3ex}

\sect{Introduction}

Let $G$ be a compact Lie group and $\g$, its Lie algebra.
Suppose $G$ acts on a symplectic manifold $(X,\om)$ in a Hamiltonian fashion,
with moment map $\mu\colon X\to\g^*$.
If $0$ is a regular value of $\mu$ and, for simplicity, if $G$ acts freely
on $\mu^{-1}(0)$, then the symplectic quotient $X_0=\mu^{-1}(0)/G$ is 
a smooth symplectic manifold.
Recall that a prequantum line bundle $L$ on $(X,\om)$ is a complex
line bundle whose curvature is $\frac{\om}{\ii}$.
Suppose such a bundle $L$ exists.
If the $G$-action lifts to $L$, then the bundle $L$ descends to 
a prequantum line bundle $L_0$ on $X_0$.

When $X$ is a \ka manifold, $L$ can be chosen as a holomorphic bundle.
If there is a holomorphic $G$-action on $X$ that lifts to $L$, then the 
twisted Dolbeault cohomology groups $H^k(X,\OO(L))$ are representations
of $G$.
The space of holomorphic sections $H^0(X,\OO(L))$ is usually considered as
the quantization of the \ka manifold $X$.
Guillemin and Sternberg [\ref{GS}] showed that if $X$ is compact and
K\"ahler, then $H^0(X,\OO(L))^G\cong H^0(X_0,\OO(L_0))$.
When $(X,\om)$ is symplectic, the individual cohomology groups do not make
sense, but there is a Spin$^c$-Dirac operator $\not\kern-5pt D$ (depending on 
a compatible almost complex structure and twisted by $L$) whose index space
$H(X)={\mathrm ker}\not\kern-4pt D\ominus{\mathrm coker}\not\kern-4pt D$
is equal to $\oplus_{k\ge 0}(-1)^k H^k(X,\OO(L))$ when $X$ is K\"ahler.
Thus $H(X)$ is a natural definition of quantization in the general situation.
The quantization conjecture of Guillemin and Sternberg is that
	\be\label{gs}
H(X)^G\cong H(X_0).
	\ee
When $X$ is compact and when $G$ is a torus or non-Abelian group,
the conjecture was proved by Meinrenken [\ref{Me1}, \ref{Me2}],
Jeffrey and Kirwan [\ref{JK}], Vergne [\ref{V}] and others under 
various generalities using localization techniques (see [\ref{Sj}] for 
a review), and by Tian and Zhang [\ref{TZ}] using an analytic approach.

If in addition $X$ is a compact \ka manifold, Tian and Zhang [\ref{TZ}] proved 
Morse-type inequalities
	\be\label{tz}
\sum_{k\ge0}\dim H^k(X_0,\OO(L_0))=\sum_{k\ge0}\dim H^k(X,\OO(L))^G+(1+t)Q(t)
	\ee
for some $Q(t)\ge0$.
Recently, it was shown by Teleman [\ref{T}], and by Braverman [\ref{Bra2}]
using a different method, that $H^k(X,\OO(L))^G\cong H^k(X_0,\OO(L_0))$
for every $k\ge0$.
Therefore the inequalities (\ref{tz}) are in fact equalities.

In this paper, we show that $H^k(X,\OO(L))^G\cong H^k(X_0,\OO(L_0))$
holds for a class of non-compact \ka manifolds with torus actions
satisfying the assumption of [\ref{PW}].
Our approach is applying the analog of instanton complex in holomorphic
Morse theory [\ref{Wu2}] to the symplectic cuts of Lerman [\ref{Le}].
This can be compared with [\ref{DGMW}], which applies the fixed-point formula
to the same geometric setting for compact manifolds.
As another application of our method, 
we settle a conjecture in [\ref{WZ}, Remark~4.11] on the exact sequence
of higher cohomology groups in the context of symplectic cutting.

For non-compact non-\ka symplectic manifolds, the validity of 
(\ref{gs}) remain open.

The rest of this paper is organized as follows.
In section~2, we recall some notations and review the main results
of holomorphic Morse theory.
In section~3, we apply the above results to the prequantum line bundles
of symplectic cuts and quotients.

\sect{Review of the instanton complex in holomorphic Morse theory}

The analog of instanton complex in holomorphic Morse theory [\ref{Wu2}]
is a spectral sequence converging to the Dolbeault cohomology groups
twisted by a holomorphic vector bundle.
This construction is possible only when the connected components of the
fixed-point set form a partially ordered set.
The existence of such a spectral sequence implies holomorphic Morse 
inequalities and fixed-point formulas on a possibly non-compact manifold.
When the manifold is compact and K\"ahler, the partial order condition 
is automatically satisfied.
Thus the results of [\ref{W2}, \ref{MW}, \ref{Wu}, \ref{WZ}] are recovered.
On the other hand, the example in [\ref{Wu}] shows that without the partial
order condition, the holomorphic Morse inequalities can fail.
Here we consider a class of non-compact \ka manifolds where the partial order
condition holds.

We first recall some notations and the basic facts of holomorphic torus 
actions.
For more details, see [\ref{Wu2}, \S 2] and references therein.

Let $T$ be a complex torus with Lie algebra $\gt$.
Let $T_\re$ be the (real) maximal compact torus subgroup of $T$
and $\gt_\re=\ii\,{\rm Lie}(T_\re)$.
Let $\ell$ be the integral lattice in $\gt_\re$, and $\ell^*\subset\gt^*_\re$,
the dual lattice.
For a (complex) representation $W$ of $T$, we let $W_\lam$ be the subspace of 
weight $\lam\in\ell^*$ and 
$\ch W=\sum_{\lam\in\ell^*}\dim_\co W_\lam\e{\lam}$, the character of $W$.
If $T=\CM$, the multiplicative group of non-zero complex numbers,
then $T_\re=S^1$, $\gt_\re=\RE$, and $\ell=\ZZ$.
In general, for any $v\in\ell-\zero$, there is an embedding $j_v\colon\CM\to T$
whose image $\CM_v$ is a $\CM$-subgroup of $T$.

Let $X$ be a complex manifold of dimension $n$.
Suppose $T$ acts holomorphically and effectively on $X$.
The fixed-point set $X^T$ of $T$ in $X$, if non-empty,
is a complex submanifold of $X$.
Let $F$ be the set of connected components of $X^T$.
Then $X^T=\bigcup_{\al\in F}X^T_\al$, where $X^T_\al$ is the component
labeled by $\al\in F$.
Let $n_\al=\dim_\co X^T_\al$.
Let $N_\al\to X^T_\al$ be the (holomorphic) normal bundle of $X^T_\al$ in $X$.
$T$ acts on $N_\al$ preserving the base $X^T_\al$ pointwise.
The weights of the isotropy representation on the normal fiber
remain constant within any connected component.
Let $\lam_{\al,k}\in\ell^*-\zero\subset\gt_\re^*$ ($1\le k\le n-n_\al$)
be the isotropy weights on $N_\al$.
The hyperplanes $(\lam_{\al,k})^\perp\subset\gt_\re$ cut $\gt_\re$ into open
polyhedral cones called {\em action chambers} [\ref{PW}].
Choose an action chamber $C$.
Let $\lam^C_{\al,k}=\pm\lam_{\al,k}$,
with the sign chosen so that $\lam^C_{\al,k}\in C^*$.
(Here $C^*$ is the dual cone in $\gt_\re^*$ defined by
$C^*=\set{\xi\in\gt_\re^*}{\bra\xi,C\ket>0}$.)
We define $\nu^C_\al$ as the number of weights $\lam_{\al,k}\in C^*$.
Let $N^C_\al$ be the direct sum of the sub-bundles corresponding to
the weights $\lam_{\al,k}\in C^*$.
Then $N_\al=N^C_\al\oplus N^{-C}_\al$.
$\nu^C_\al$ is the rank of the holomorphic vector bundle $N^C_\al$;
that of $N^{-C}_\al$ is $\nu^{-C}_\al=n-n_\al-\nu^C_\al$,
which is called the {\em polarizing index} of $X^T_\al$ with respect to $C$.

We further assume that $(X,\om)$ is a \ka manifold and that the action of
the compact torus $T_\re$ is Hamiltonian with respect to the \ka form $\om$.
Then $X^T_\al$ ($\al\in F$) are \ka submanifolds.
Let $\mu\colon X\to\gt^*_\re$ be the moment map.

\begin{assump}\label{PWA}
{\em ([\ref{PW}, Assumption~1.3])
There is $v\in\gt_\re$ such that $\bra\mu,v\ket\colon X\to\RE$ is proper 
and not surjective and $F$ is a (non-empty) finite set.}
\end{assump}

Under this assumption, there is an action chamber $C$ such that the function
$\bra\mu,v\ket$ is bounded from above if $v\in C$ [\ref{PW}].
Since for any $x\in X$, $\bra\mu,v\ket$ strictly decreases under the flow
$t\mapsto j_v(\e{t})x$, the limit $\lim_{u\to0}j_v(u)x$ exists.
This limit depends only on the action chamber $C$ containing $v$ and
is therefore denoted by $\pi^C(x)$.
Clearly $\pi^C\colon X\to X^T$.
Set $X^C_\al=(\pi^C)^{-1}(X^T_\al)$.
Then we have the {\em \bb decomposition} $X=\bigcup_{\al\in F}X^C_\al$.
If for $x\in X$, the limit $\lim_{u\to\infty}j_v(u)x$ also exists, 
it depends only on $C$ as well; we denote it by $\pi^{-C}(x)$.

We define a relation $\prec$ on $F$.
For $\al,\beta\in F$, we write $\al\to\beta$ if there is $x\in X$
such that $\pi^C(x)\in X^T_\al$ and $\pi^{-C}(x)\in X^T_\beta$.
We write $\al\prec\beta$ if either $\al=\beta$ or there is a {\em chain} 
from $\al$ to $\beta$,
i.e., a finite sequence $\al_0=\al, a_1,\dots,\al_{r-1},\al_r=\beta$
in $F$ such that $\al_{i-1}\to\al_i$ for all $1\le i\le r$ ($r>0$).
For general holomorphic torus actions, $(F,\prec)$ need not be a partially
ordered set. 
Holomorphic Morse theory fails precisely in this case [\ref{Wu}, \ref{Wu2}].
In our situation of non-compact \ka manifolds, $\prec$ is a partial ordering
since $\bra\mu,v\ket$ is a strictly decreasing function on $(F,\prec)$.

A consequence of having a partial ordering $\prec$ on $F$ is as follows.
The \bb decomposition becomes {\em filterable} in the sense that there is
a descending sequence of $T$-invariant subvarieties
$X=Z_0\supset Z_1\supset\cdots\supset Z_m\supset Z_{m+1}=\emptyset$
such that for all $0\le p\le m$, $Z_p-Z_{p+1}=\bigcup_{\al\in F_p}X^C_\al$
for a subset $F_p\subset F$ such that neither $\al\prec\beta$ nor
$\beta\prec\al$ if $\al\ne\beta\in F_p$.

Now we state the main result of [\ref{Wu2}] specialized to \ka manifolds
satisfying Assumption~\ref{PWA}.
The following is a combination of [\ref{Wu2}, Proposition~2.16] and
[\ref{Wu2}, Theorem~3.14.1].

\begin{thm}\label{MAIN'}
Let $X$ be a \ka manifold with an effective holomorphic $T$-action such that
the $T_\re$-action is Hamiltonian and satisfies Assumption~\ref{PWA}.
Let $E$ be a holomorphic vector bundle over $X$ on which the $T$-action lifts
holomorphically.
Then there is a $T$-equivariant spectral sequence converging $T$-equivariantly
to $H^*(X,\OO(E))$ with
	\be\label{e1'}
E\pq_1=
\bigoplus_{\al\in F_p}H^{p+q+\nu^C_\al+n_\al-n}(X^T_\al,\OO(S((N^C_\al)^*)
\otimes S(N^{-C}_\al)\otimes\medwedge^{n-n_\al-\nu^C_\al}(N^{-C}_\al)
\otimes E|_{X^T_\al})).
	\ee
\end{thm}

We remark that there is another spectral sequence converging to the
Dolbeault cohomology groups with compact support [\ref{Wu2}, Theorem~3.7.1],
which seems not to be related to problems of geometric quantization
[\ref{Wu2}, Remark~4.9.1].

\sect{Applications to higher cohomology groups of \ka quotients}

We apply Theorem~\ref{MAIN'} to the prequantum line bundles on \ka manifolds.
Let $(X,\om)$ be a (possibly non-compact) \ka manifold with a holomorphic
action of the complex torus $T$.
Assume that the action of the compact torus $T_\re$ is Hamiltonian and 
the moment map $\mu\colon X\to\gt_\re^*$ satisfies Assumption~\ref{PWA}.
Recall that a prequantum line bundle $L$ on $(X,\om)$ is a holomorphic
line bundle whose curvature is $\frac{\om}{\ii}$.
Suppose such a bundle $L$ exists and the $T$-action lifts to a holomorphic 
action on $L$.
We choose the lifted action such that the weight on the fiber $L|_{X^T_\al}$
($\al\in F$) is $\mu(X^T_\al)\in\ell^*$.
We will study the cohomology groups $H^k(X,\OO(L))$ ($k\ge0$) as
representations of $T$.
Suppose $0$ is a regular value of $\mu$.
For simplicity, we assume that the $T_\re$-action on $\mu^{-1}(0)$ is free.
Then the symplectic quotient $X_0=\mu^{-1}(0)/T_\re$ is a smooth \ka manifold.
Moreover, $L$ descends to a prequantum line bundle $L_0$ on $X_0$.

We first consider the case $T=\CM$.
Without loss of generality, we assume that the moment map $\mu$ 
is bounded from above.
Choose the action chamber $C=\RE^+$ and set $N^\pm_\al=N^{\pm C}_\al$,
$\nu_\al=\nu^C_\al$ ($\al\in F$).
We construct the symplectic cuts $(X_\pm,\om_\pm)$ as the symplectic quotients
of the $S^1$-action on $X\times\CO$, where the weights on $\CO$ are $\pm1$,
respectively [\ref{Le}].
The two cuts $X_\pm$ are \ka manifolds with holomorphic $\CM$-actions.
$X_+$ is compact and $X_-$ satisfies Assumption~\ref{PWA}.
The sets of connected components of $X^\cm_\pm$ are
$F_\pm=\zero\cup\set{\al\in F}{\mu(X^\cm_\al)\in\RE^\pm}$, respectively,
and $X^\cm_{\pm,0}\cong X_0$, $X^\cm_{\pm,\al}\cong X^\cm_\al$
as complex manifolds [\ref{WZ}, Lemma~4.6], which we now identify.
Let $N_0\to X_0$ be the holomorphic line bundle associate to the
circle bundle $\mu^{-1}(0)\to X_0$.
Then $\CM$ acts on the fibers of $N_0$ with weight $1$.
The holomorphic normal bundles of $X_0$ in $X_\pm$ are isomorphic to
$N_0^{\mp1}$, respectively.
Since the action of $\CM$ lifts to $L$, the prequantum line bundles
$L_0\to X_0$ and $L_\pm\to X_\pm$ exist.
We have the isomorphisms $L_\pm|_{X_0}\cong L_0$ and
$L_\pm|_{X_\pm-X_0}\cong L|_{\mu^{-1}(\re^\pm)}$
(see for example [\ref{WZ}, Lemma~4.9]).

\begin{thm}\label{CUTS} Under the above assumptions, we have, 
for every $k\ge 0$,\\
1.
	\be
H^k(X_\pm,\OO(L_\pm))^\cm\cong H^k(X_0,\OO(L_0)).
	\ee
2.
	\be
H^k(X_+,\OO(L_+))_\lam\cong\two{H^k(X,\OO(L))_\lam}{\lam\ge0}{0}{\lam<0.}
	\ee
3. if $S^1$ acts on $\mu^{-1}(a)$ freely for sufficiently negative $a<0$, then
	\be
H^k(X_-,\OO(L_-))_\lam\cong\two{0}{\lam>0}{H^k(X,\OO(L))_\lam}{\lam\le0.}
	\ee
\end{thm}

\proof{1. Recall that the normal bundle of $X_0$ in $X_-$ is $N_0$.
By Theorem~\ref{MAIN'}, $H^k(X_-,\OO(L_-))$ can be computed by a spectral
sequence whose $E_1$-terms are given by the cohomology groups
$H^*(X_\al^\cm,\OO(S((N^+_\al)^*)\otimes S(N^-_\al)\otimes
\medwedge^{n-n_\al-\nu_\al}(N^-_\al)\otimes L|_{X^{\smco^\times}_\al}))$
($\al\in F_--\zero$) and $H^*(X_0,\OO(S(N_0^{-1})\otimes L_0))$.
Since the weight of $\CM$ in $L|_{X^{\smco^\times}_\al}$ is 
$\mu(X^\cm_\al)<0$ for $\al\in F_--\zero$,
the weights in the former cohomology groups (on $X^\cm_\al$) are negative.
On the other hand, the weights in the latter (on $X_0$) are non-positive;
the $\CM$-invariant part (with weight $0$) comes from $H^*(X_0,\OO(L_0))$.
Therefore, $H^k(X_-,\OO(L_-))_\lam=0$ for $\lam>0$.
Moreover, the $\CM$-invariant part of the spectral sequence degenerates
at $E_1$, and hence $H^k(X_-,\OO(L_-))^\cm\cong H^k(X_0,\OO(L_0))$.
By reversing the $S^1$-action, we obtain $H^k(X_+,\OO(L_+))_\lam=0$ for 
$\lam<0$ and $H^k(X_+,\OO(L_+))^\cm\cong H^k(X_0,\OO(L_0))$.\\
2. Next, we show that $H^k(X_+,\OO(L_+))_\lam\cong H^k(X,\OO(L))_\lam$
for $\lam\ge0$.
The normal bundle of $X_0$ in $X_+$ is $N_0^{-1}$.
By Theorem~\ref{MAIN'},
the $E_1$-terms of the spectral sequence for $H^k(X_+,\OO(L_+))$ are given by
$H^*(X_\al^\cm,\OO(S((N^+_\al)^*)\otimes S(N^-_\al)\otimes
\medwedge^{n-n_\al-\nu_\al}(N^-_\al)\otimes L|_{X^{\smco^\times}_\al}))$
(with $\al\in F_+-\zero$) and 
$H^*(X_0,\OO(S(N_0^{-1})\otimes N_0^{-1}\otimes L_0))$,
whereas those for $H^k(X,\OO(L))$ are given by
$H^*(X_\al^\cm,\OO(S((N^+_\al)^*)\otimes S(N^-_\al)\otimes
\medwedge^{n-n_\al-\nu_\al}(N^-_\al)\otimes L|_{X^{\smco^\times}_\al}))$
with all $\al\in F$.
For both spectral sequences, the restrictions of the $E_1$-terms to a
non-negative weight $\lam\ge0$ are $H^*(X_\al^\cm,\OO(S((N^+_\al)^*)
\otimes S(N^-_\al)\otimes\medwedge^{n-n_\al-\nu_\al}(N^-_\al)
\otimes L|_{X^{\smco^\times}_\al}))_\lam$ ($\al\in F_+-\zero$).
{}From the proof of [\ref{Wu2}, Proposition~2.16], 
the open sets $X_+-X_0\subset X_+$ and $\set{x\in X}{\mu(x)>0}\subset X$
are holomorphically isomorphic.
So for the two spectral sequences associated to $X_+$ and $X$, the coboundary 
operators are the same after the restriction to the non-negative weight spaces.
This implies the desired isomorphism.\\
3. Finally, we show that $H^k(X_-,\OO(L_-))_\lam\cong H^k(X,\OO(L))_\lam$
for $\lam\le0$.
If $X$ is compact, this follows from part~2 by reversing the $S^1$-action.
If $\mu$ is not bounded from below, choose a sufficiently negative $a<0$.
We shift the moment map to $\mu-a$ and denote by $X_{\ge a}$ the symplectic cut
above level $a$ and by $L_{\ge a}$ the prequantum line bundle over $X_{\ge a}$.
Under our assumptions, $X_{\ge a}$ and $(X_-)_{\ge a}$ are smooth compact \ka 
manifolds.
Using part~2, we conclude that 
$H^k(X,\OO(L))_\lam\cong H^k(X_{\ge a},\OO(L_{\ge a}))_\lam$ and
$H^k(X_-,\OO(L_-))_\lam\cong H^k((X_-)_{\ge a},\OO((L_-)_{\ge a}))_\lam$
whenever $\lam>a$.
The proof is thus reduced to the compact situation.}

\begin{rmk}{\em
1. Theorem~\ref{CUTS}.1 implies that the Morse-type inequalities 
in [\ref{WZ}, Proposition~4.10] are equalities.\\
2. We believe that the further assumption in Theorem~\ref{CUTS}.3 
in the non-compact situation is technical can be avoided 
in a more careful analysis.}
\end{rmk}

We now consider the case of torus action.
Recall the notations at the beginning of this section.

\begin{cor}\label{HIGHER}
Let $X$ be a (possibly non-compact) \ka manifold with a holomorphic $T$-action.
Suppose the $T_\re$-action is Hamiltonian and satisfies Assumption~\ref{PWA}.
Suppose $0$ is a regular value of the moment map $\mu$ and $T_\re$ acts on
$\mu^{-1}(0)$ freely.
Then there is an isomorphism
	\be\label{higher}
H^k(X,\OO(L))^T\cong H^k(X_0,\OO(L_0))
	\ee
for every $k\ge0$.
\end{cor}

\proof{From parts 1 and 2 of Theorem~\ref{CUTS}, we get
$H^k(X,\OO(L))^\cm\cong H^k(X_+,\OO(L_+))^\cm\cong H^k(X_0,\OO(L_0))$;
this is (\ref{higher}) when $T=\CM$.
The general case follows by induction using reduction in stages.}

\begin{ex}{\em
We give an example in which a higher cohomology group does not vanish.
Let $\Sig$ be a Riemann surface of genus $g\ge2$ with constant curvature.
Then the holomorphic cotangent bundle $T^*\Sig$ is the prequantum line bundle
for a symplectic form on $\Sig$.
We have $\dim_\co H^0(\Sig,\OO(T^*\Sig))=g$ and 
$\dim_\co H^1(\Sig,\OO(T^*\Sig))=1$.
Consider $\CP^1$ with the standard $\CM$-action such that the image of 
the moment map is $[-1,1]$.
Let $L_1$ be prequantum line bundle over $\CP^1$.
Let $X=\Sig\times\CP^1$ so that $\CM$ acts only on the second factor.
Let $\pi_\Sig$ and $\pi_1$ be the projections from $X$ to $\Sig$ and $\CP^1$,
respectively.
Then $L=\pi_\Sig^* T\Sig\otimes\pi_1^*L_1$ is a prequantum line bundle over 
$X$.
The action of $S^1$ on $X$ is Hamiltonian and the symplectic quotient
$X_0\cong\Sig$, with $L_0\cong T^*\Sig$.
According to Corollary~\ref{HIGHER},
$H^1(X,\OO(L))^\cm\cong H^1(\Sig,\OO(T^*\Sig))\neq0$.
For an example in which the manifold is not a product, consider
the projectivization of a non-trivial line bundle over $\Sig$.}
\end{ex}

Finally, we return to the study of symplectic cuts of $\CM$-actions.

\begin{cor}\label{CONJ}
Under the conditions of Theorem~\ref{CUTS},
there is an $\CM$-equivariant short exact sequence
	\be\label{short}
0\to H^k(X,\OO(L))\to H^k(X_+,\OO(L_+))\oplus H^k(X_-,\OO(L_-))\to
H^k(X_0,\OO(L_0))\to 0
	\ee
for every $k\ge0$.
\end{cor}

\proof{We decompose (\ref{short}) into weight spaces of the $\CM$-action.
For a positive or negative weight $\lam\in\ZZ$, (\ref{short}) reduces to
$H^k(X,\OO(L))_\lam\cong H^k(X_+,\OO(L_+))_\lam$ and
$H^k(X,\OO(L))_\lam\cong H^k(X_-,\OO(L_-))_\lam$, respectively.
The $\CM$-invariant part of (\ref{short}) is exact because of the isomorphisms
$H^k(X,\OO(L))^\cm\cong H^k(X_\pm,\OO(L_\pm))^\cm\cong H^k(X_0,\OO(L_0))$.}

\begin{rmk}{\em
In [\ref{TZ}, Remark~4.11], it was conjectured that when $M$ is compact,
there is an $\CM$-equivariant long exact sequence
	\be\label{WZseq}
\cdots\to H^k(X,\OO(L))\to H^k(X_+,\OO(L_+))\oplus H^k(X_-,\OO(L_-))\to
H^k(X_0,\OO(L_0))\to H^{k+1}(X,\OO(L))\to\cdots
	\ee
This implies, among other things, Morse-type inequalities
	\bea\label{WZineq}
\vc\sum_{k=0}^nt^k\ch H^k(X_+,\OO(L_+))+\sum_{k=0}^nt^k\ch H^k(X_-,\OO(L_-))
   +\sum_{k=0}^{n-1} t^{k+1}\ch H^k(X_0,\OO(L_0))                       \nno
\eq\quad\sum_{k=0}^nt^k\ch H^k(X,\OO(L))+(1+t)Q(t),
        \eea
for some character-valued polynomial $Q(t)\ge0$.
Braverman [\ref{Bra1}] proved that when $M$ is compact, (\ref{WZineq})
holds for a general vector bundle on which the $\CM$-action lifts;
the proof was based on an elegant construction of a holomorphic family
of complex manifolds which degenerate at one fiber to the union of $X_+$
and $X_-$ along $X_0$.
However the long exact sequence (\ref{WZseq}) for the prequantum line bundle
has remained open except in some examples where all the cohomology groups are
explicitly known.
(See for example [\ref{Bra1}, Example~2.14].)
Corollary~\ref{CONJ} is stronger than exactness of (\ref{WZseq}) and holds
for possibly non-compact manifolds.}
\end{rmk}

\bigskip
\newpage
        \newcommand{\athr}[2]{{#1}.\ {#2}}
        \newcommand{\au}[2]{\athr{{#1}}{{#2}},}
        \newcommand{\an}[2]{\athr{{#1}}{{#2}} and}
        \newcommand{\jr}[6]{{#1}, {\it {#2}} {#3}\ ({#4}) {#5}-{#6}}
        \newcommand{\pr}[3]{{#1}, {#2} ({#3})}
        \newcommand{\bk}[4]{{\it {#1}}, {#2}, ({#3}, {#4})}
        \newcommand{\cf}[8]{{\it {#1}}, {#2}, {#5},
                 {#6}, ({#7}, {#8}), pp.\ {#3}-{#4}}
        \vspace{5ex}
        \begin{flushleft}
{\bf References}
        \end{flushleft}
{\small
        \begin{enumerate}

	\item\label{Bra1}
	\au{M}{Braverman}
	\pr{Symplectic cutting of \ka manifolds}
	{preprint {\tt alg-geom/9712024}}{December, 1997}

	\item\label{Bra2}
	\au{M}{Braverman}
	\pr{Cohomology of the Mumford quotients}
	{preprint {\tt math.SG/9809146}}{September, 1998}
	
        \item\label{DGMW}
        \au{H}{Duistermaat} \au{V}{Guillemin} \an{E}{Meinrenken} \au{S}{Wu}
        \jr{Symplectic reduction and Riemann-Roch for circle actions}
        {Math.\ Res.\ Lett.}{2}{1995}{259}{266}

        \item\label{GS}
        \an{V}{Guillemin} \au{S}{Sternberg}
        \jr{Geometric quantization and multiplicities of group representations}
        {Invent. Math.}{67}{1982}{515}{538}

        \item\label{JK}
        \an{L.\ C}{Jeffrey} \au{F.\ C}{Kirwan}
        \jr{Localization and the quantization conjecture}
        {Topology}{36}{1997}{647}{693}

        \item\label{Le}
        \au{E}{Lerman}
        \jr{Symplectic cuts}{Math.\ Res.\ Lett.}{2}{1995}{247}{258}

	\item\label{MW}
	\an{V}{Mathai} \au{S}{Wu}
	\jr{Equivariant holomorphic Morse inequalities I: a heat kernel proof}
	{J.\ Diff.\ Geom.}{46}{1997}{78}{98}

        \item\label{Me1}
        \au{E}{Meinrenken}
        \jr{On Riemann-Roch formulas for multiplicities}
        {J.\ Amer.\ Math.\ Soc.}{9}{1996}{373}{389}

        \item\label{Me2}
        \au{E}{Meinrenken}
        \jr{Symplectic surgery and the Spin$^c$-Dirac operator}
        {Adv.\ Math.}{134}{1998}{240}{277}

        \item\label{PW}
        \an{E}{Prato} \au{S}{Wu}
        \jr{Duistermaat-Heckman measures in a non-compact setting}
        {Comp.\ Math.}{94}{1994}{113}{128}

	\item\label{Sj}
	\au{R}{Sjamaar}
	\jr{Symplectic reduction and Riemann-Roch formulas for multiplicities}
	{Bull.\ Amer.\ Math.\ Soc.}{33}{1996}{327}{338}

	\item\label{T}
	\au{C}{Teleman}
	\pr{The quantization conjecture revisited}
	{preprint {\tt math.AG/9808029}}{August, 1998}

        \item\label{TZ}
        \an{Y}{Tian} \au{W}{Zhang}
        \jr{Symplectic reduction and quantization}
        {C. R. Acad. Sci. Paris, S\'erie I}{324}{1997}{433}{438};
        \jr{An analytic proof of the geometric quantization conjecture of
        Guillemin-Sternberg}{Invent.\ Math.}{132}{1998}{229}{259}

        \item\label{V}
        \au{M}{Vergne}
        \jr{Multiplicities formula for geometric quantization. I, II}
        {Duke Math. J.}{82}{1996}{143}{179}, 181-194

	\item\label{W2} 
	\au{E}{Witten} 
	\cf{Holomorphic Morse inequalities}
	{Algebraic and differential topology, Teubner-Texte Math., 70}
	{318}{333}{ed.\ G.\ Rassias}{Teubner}{Leipzig}{1984}

	\item\label{Wu}
	\au{S}{Wu}
	\pr{Equivariant holomorphic Morse inequalities II: torus and
	non-Abelian group actions}
	{MSRI preprint No.~1996-013, {\tt dg-ga/9602008}}{1996}

	\item\label{Wu2}
	\au{S}{Wu}
	\pr{On the instanton complex of holomorphic Morse theory}
	{Adelaide IGA preprint 1998-07, {\tt math.AG/9806168}}{June, 1998}

        \item\label{WZ}
        \an{S}{Wu} \au{W}{Zhang}
        \jr{Equivariant holomorphic Morse inequalities III: non-isolated
        fixed points}{Geom.\ Funct.\ Anal.}{8}{1998}{149}{178}

        \end{enumerate}}
        \end{document}